\documentclass[12pt]{article}
\usepackage{amssymb}
\usepackage{amsmath}

\input{tcilatex}
\begin{document}

\begin{center}
\textbf{Iterated function systems consisting of }${\small \varphi }$\textbf{%
-max-contractions have attractor}

\bigskip

by \textit{Flavian GEORGESCU, Radu MICULESCU }and\textit{\ Alexandru MIHAIL}

\bigskip
\end{center}

\textbf{Abstract}. {\small We associate to each iterated function system
consisting of }${\small \varphi }${\small -max-contractions an operator (on
the space of continuous functions from the shift space on the metric space
corresponding to the system) having a unique fixed point whose image turns
out to be the attractor of the system. Moreover, we prove that the unique
fixed point of the operator associated to an iterated function system
consisting of convex contractions is the canonical projection from the shift
space on the attractor of the system.}

\bigskip



\bigskip

\textbf{1. Introduction}

\bigskip

The importance of the concepts of shift space (or the code space) and
canonical projection associated to an iterated function system in the
description of topological properties of the attractor of such a system was
pointed out in several papers like [1] (where the theory of fractal tops is
treated), [5] (where the shift space and the canonical projection for an
infinite iterated function system are studied) and [3]. A special place in
this discussion deserves the paper [8] where the canonical projection
between the shift space of an infinite iterated function system and its
attractor is presented as a fixed point, in two cases: a) the constitutive
functions of the system are uniformly Meir-Keeler; b) the metric space
associated to the system is compact and the system consists of a finite
number of contractive functions.

As part of the current effort to extend Hutchinson theory concerning
iterated function systems to more general frameworks, in [7], the concept of
iterated function system consisting of convex contractions was introduced
and the existence and uniqueness of the attractor of such a system was
obtained. See also [4] for a more general result in this direction.

In this paper we introduce the concept of iterated function system
consisting of $\varphi $-$\max $-contractions -for short $\varphi $-$\max $%
-IFS- (see Definition 2.12) which generalizes the concept of iterated
function system consisting of convex contractions (see Definition 2.13). To
such a system $\mathcal{S}$ we associate an operator $G_{\mathcal{S}}:%
\mathcal{C}\rightarrow \mathcal{C}$, where $\mathcal{C}$ stands for the
space of continuous functions from the shift space on the metric space
corresponding to the system, which has a unique fixed point (see Theorem
3.1) whose image is the attractor of $\mathcal{S}$ (see Theorem 3.2). In
this way we provide a new method to prove the existence and uniqueness of
the attractor of an iterated function system since the classical approach
consists on proving that $F_{\mathcal{S}}$ -the fractal operator associated
to the system $\mathcal{S}$ whose constitutive functions belong to a fixed
family $\mathcal{F}$ of Picard contractions- is also an element of $\mathcal{%
F}$.

Moreover, we point out the following two facts:

\medskip

$\bullet $ In the particular case of an iterated function system consisting
of convex contractions the fixed point of $G_{\mathcal{S}}$ turns out to be
the canonical projection from the shift space on the attractor of the system
(see Theorem 3.4). This is a companion result of those from [8].

\medskip

$\bullet $ Since Theorem 3.1 is also valid for a class of iterated function
systems involving a possible infinite family of $\varphi $-$\max $%
-contractions, abbreviated $\varphi $-$\max $-PIIFSs (see Definition 4.1 and
Theorem 4.4), we raise the following open question: Is it true that the
fractal operator associated to a $\varphi $-$\max $-PIIFS $\mathcal{S}$ is a
Picard operator whose fixed point is the image of the fixed point of the
operator $G_{\mathcal{S}}$?

\bigskip

\textbf{2.} \textbf{Preliminaries}

\bigskip

\textbf{Some notations}

\bigskip

Given the sets $A$ and $B$, by $B^{A}$ we mean the set of functions from $A$
to $B$.

\bigskip

Given a function $f:X\rightarrow X$ and $p\in \mathbb{N}$, by $f^{[p]}$ we
mean $\underset{p\text{ times}}{f\circ f\circ ...\circ f}$.

\bigskip

Given a metric space $(X,d)$, by:

- $P_{cp}(X)$ we mean the set of non-empty compact subsets of $X$

- $P_{b,cl}(X)$ we mean the set of bounded and closed subsets of $X$.

\bigskip

Given a metric space $(X,d)$, a subset $A$ of $X$ and $\varepsilon >0$, by $%
E_{\varepsilon }(A)$, the $\varepsilon $-expansion of $A$, we mean 
\begin{equation*}
\{y\in X\mid \text{there exists}\mathit{\ }x\in A\mathit{\ }\text{such that\ 
}d(x,y)<\varepsilon \}\text{.}
\end{equation*}

\newpage

\textbf{The Hausdorff-Pompeiu metric}

\bigskip

\textbf{Definition 2.1.} \textit{The function }$H:P_{cp}(X)\times
P_{cp}(X)\rightarrow \lbrack 0,+\infty )$\textit{, where }$(X,d)$\textit{\
is a metric space,\ defined by}

\begin{equation*}
H(A,B)=\max (d(A,B),d(B,A))=
\end{equation*}

\begin{equation*}
=\inf \{\varepsilon \in \lbrack 0,\infty )\mid A\subseteq E_{\varepsilon }(B)%
\text{ \textit{and} }B\subseteq E_{\varepsilon }(A)\}\text{,}
\end{equation*}%
\textit{where }$d(A,B)=\underset{x\in A}{\sup }d(x,B)=\underset{x\in A}{\sup 
}(\underset{y\in B}{\inf }d(x,y))$, \textit{turns out to be a metric which
is called the Hausdorff-Pompeiu metric.}

\bigskip

\textbf{Remark 2.2.}

\textit{a)} \textit{Given two metric spaces }$(X,d)$ \textit{and} $%
(Y,d^{^{\prime }})$\textit{,\ a sequence of continuous functions }$%
(f_{n})_{n\in \mathbb{N}}$\textit{, where }$f_{n}:X\rightarrow Y$\textit{,
and }$K\in P_{cp}(X)$, \textit{we have }%
\begin{equation*}
\underset{n\rightarrow \infty }{\lim }H(f_{n}(K),f(K))=0\text{\textit{,}}
\end{equation*}%
\textit{\ provided that }$f_{n}\overset{u}{\rightarrow }f$\textit{.}

Indeed, let us fixed an arbitrary $\varepsilon >0$. Since $f_{n}\overset{u}{%
\rightarrow }f$, there exists $n_{\varepsilon }\in \mathbb{N}$ such that $%
d(f_{n}(x),f(x))<\varepsilon $ for every $x\in K$ and every $n\in \mathbb{N}$%
, $n\geq n_{\varepsilon }$. Hence $f_{n}(K)\subseteq E_{\varepsilon }(f(K))$
and $f(K)\subseteq E_{\varepsilon }(f_{n}(K))$, so $H(f_{n}(K),f(K))<%
\varepsilon $ for every $n\in \mathbb{N}$, $n\geq n_{\varepsilon }$.
Consequently $\underset{n\rightarrow \infty }{\lim }H(f_{n}(K),f(K))=0$.

\textit{b) }(see Proposition 2.8 from [7])\textit{\ Given a complete metric
space} $(X,d)$, $(Y_{n})_{n\in \mathbb{N}}\subseteq P_{cp}(X)$ \textit{and }$%
Y\in P_{cp}(X)$\textit{, we have}%
\begin{equation*}
Y\cup (\underset{n\in \mathbb{N}}{\cup }Y_{n})\in P_{cp}(X)\text{,}
\end{equation*}%
\textit{provided that} $\underset{n\rightarrow \infty }{\lim }H(Y_{n},Y)=0$.

\bigskip

\textbf{Comparison functions}

\bigskip

\textbf{Definition 2.3.} \textit{A function }$\varphi :[0,\infty
)\rightarrow \lbrack 0,\infty )$\textit{\ is called a comparison function
provided that it satisfies the following properties:}

\textit{i) }$\varphi $\textit{\ is increasing;}

\textit{ii) }$\underset{n\rightarrow \infty }{\lim }\varphi ^{\lbrack
n]}(x)=0$\textit{\ for every }$x\in \lbrack 0,\infty )$\textit{.}

\bigskip

\textbf{Remark 2.4.} \textit{For each comparison function the following two
properties are valid:}

\textit{a) }$\varphi (0)=0$\textit{;}

\textit{b) }$\varphi (x)<x$\textit{\ for every }$x\in (0,\infty )$\textit{.}

\bigskip

\textbf{A fixed point result}

\bigskip

\textbf{Definition 2.5. }\textit{Given a metric space }$(X,d)$, \textit{a
function} $f:X\rightarrow X$\textit{\ is called Picard operator if there
exists a unique fixed point }$\alpha $ \textit{of }$f$\textit{\ and the
sequence }$(f^{[n]}(x))_{n\in \mathbb{N}}$\textit{\ is convergent to }$%
\alpha $ \textit{for every }$x\in X$.

\bigskip

\textbf{Theorem 2.6} (see Theorem 3.1 from [6])\textbf{.}\textit{\ Every
continuous function }$f:X\rightarrow X$\textit{, where }$(X,d)$\textit{\ is
a complete metric space, is a Picard operator provided that there exist a
comparison function }$\varphi :[0,\infty )\rightarrow \lbrack 0,\infty )$%
\textit{\ and }$p\in \mathbb{N}^{\ast }$ \textit{such that}

\begin{equation*}
d(f^{[p]}(x),f^{[p]}(y))\leq \varphi (\underset{j\in \{0,1,2,...,p-1\}}{\max 
}d(f^{[j]}(x),f^{[j]}(y)))\text{,}
\end{equation*}%
\textit{for every }$x,y\in X$\textit{.}

\bigskip

\textbf{The shift space}

\bigskip

Given a nonempty set $I$, we denote the set $I^{\mathbb{N}^{\ast }}$ by $%
\Lambda (I)$. Thus $\Lambda (I)$ is the set of infinite words with letters
from the alphabet $I$ and a standard element $\omega $ of $\Lambda (I)$ can
be presented as $\omega =\omega _{1}\omega _{2}...\omega _{n}\omega
_{n+1}... $ .

We endow $\Lambda (I)$ with the metric described by%
\begin{equation*}
d_{\Lambda }(\omega ,\theta )=\{%
\begin{array}{cc}
0\text{,} & \text{if }\omega =\theta \\ 
\frac{1}{2^{m}}\text{,} & \text{if }\omega \neq \theta%
\end{array}%
\end{equation*}%
where, if $\omega =\omega _{1}\omega _{2}\omega _{3}...\omega _{n}\omega
_{n+1}...\neq \theta =\theta _{1}\theta _{2}\theta _{3}...\theta _{n}\theta
_{n+1}...$, $m$ is the unique natural number such that $\omega _{1}=\theta
_{1}$, $\omega _{2}=\theta _{2}$, ..., $\omega _{m-1}=\theta _{m-1}$ and $%
\omega _{m}\neq \theta _{m}$.

\bigskip

\textbf{Remark 2.7.}

\textit{a) The convergence in the metric space }$(\Lambda (I),d_{\Lambda })$%
\textit{\ is the convergence on components.}

\textit{b) }$(\Lambda (I),d_{\Lambda })$\textit{\ is a complete metric space.%
}

\textit{c) If }$I$\textit{\ is finite, then} $(\Lambda (I),d_{\Lambda })$%
\textit{\ is compact.}

\bigskip

\textbf{Proposition 2.8} (see [2], Exercise 3.6.2, page 46)\textbf{.} 
\textit{If }$I$\textit{\ is a finite set having at least two elements, then }%
$(\Lambda (I),d_{\Lambda })$\textit{\ is a Cantor set (i.e. it is compact,
perfect and totally disconnected).}

\bigskip

\textbf{Proposition 2.9} (see [10], Theorem 30.3, page 216)\textbf{.} 
\textit{Any two Cantor sets are homeomorphic.}

\bigskip

\textbf{Proposition 2.10} (see [10], Theorem 30.7, page 217)\textbf{.} 
\textit{Any compact metric space is a continuous image of a Cantor set.}

\bigskip

Taking into account the above three results we obtain the following:

\bigskip

\textbf{Proposition 2.11.} \textit{Given a finite set }$I$\textit{\ having
at least two elements and a metric space }$(X,d)$\textit{, for every }$K\in
P_{cp}(X)$ \textit{there exists a continuous function }$g_{K}:\Lambda
(I)\rightarrow X$\textit{\ such that }$g_{K}(\Lambda (I))=K$\textit{.}

\bigskip

\textbf{More notations}

\bigskip

Given a nonempty set $I$, we denote the set $I^{\{1,2,...,n\}}$ by $\Lambda
_{n}(I)$. Thus $\Lambda _{n}(I)$ is the set of words of length $n$ with
letters from the alphabet $I$ and a standard element $\omega $ of $\Lambda
(I)$ can be presented as $\omega =\omega _{1}\omega _{2}...\omega _{n}$. By $%
\Lambda _{0}(I)$ we mean the set having only one element, namely the empty
word denoted by $\lambda $.

\bigskip

For $n\in \mathbb{N}^{\ast }$, we denote by $V_{n}(I)$ the set $\underset{%
k\in \{0,1,2,...,n-1\}}{\cup }\Lambda _{k}(I)$.

\bigskip

Given a nonempty set $I$, $m,n\in \mathbb{N}$ and two words $\omega =\omega
_{1}\omega _{2}...\omega _{n}\in \Lambda _{n}(I)$\ and $\theta =\theta
_{1}\theta _{2}...\theta _{m}\in \Lambda _{m}(I)$ or $\theta =\theta
_{1}\theta _{2}...\theta _{m}\theta _{m+1}...\in \Lambda (I)$, by $\omega
\theta $ we mean the concatenation of the words $\omega $ and $\theta $, i.e.%
$\ \omega \theta =\omega _{1}\omega _{2}...\omega _{n}\theta _{1}\theta
_{2}...\theta _{m}$ and respectively $\omega \theta =\omega _{1}\omega
_{2}...\omega _{n}\theta _{1}\theta _{2}...\theta _{m}\theta _{m+1}...$ .

\bigskip

For a family of functions $(f_{i})_{i\in I}$, where $f_{i}:X\rightarrow X$,
and $\omega _{1},\omega _{2},...,\omega _{n}\in I$, by $f_{\omega _{1}\omega
_{2}...\omega _{n}}$ we mean $f_{\omega _{1}}\circ f_{\omega _{2}}\circ
...\circ f_{\omega _{n}}$. For a function $f:X\rightarrow X$, by $f_{\lambda
}$ we mean $Id_{X}$.

\bigskip

\textbf{Iterated function systems consisting of }$\varphi $-\textbf{%
max-contractions}

\bigskip

\textbf{Definition 2.12.\ }\textit{An iterated function system consisting of 
}$\varphi $-$\max $\textit{-contractions} \textit{(}$\varphi $-$\max $-%
\textit{IFS for short) is described by:}

-\textit{\ a complete metric space }$(X,d)$

-\textit{\ a finite family of continuous functions }$(f_{i})_{i\in I}$%
\textit{, where }$f_{i}:X\rightarrow X$,\textit{\ having the property that
there exist a comparison function }$\varphi :[0,\infty )\rightarrow \lbrack
0,\infty )$\textit{\ and }$p\in \mathbb{N}^{\ast }$ \textit{such that}%
\begin{equation*}
\underset{\omega \in \Lambda _{p}(I)}{\max }d(f_{\omega }(x),f_{\omega
}(y))\leq \varphi (\underset{\omega \in V_{p}(I)}{\max }d(f_{\omega
}(x),f_{\omega }(y)))\text{,}
\end{equation*}%
\textit{for every }$x,y\in X$\textit{.}

\medskip

\textit{We denote such a system by }%
\begin{equation*}
\mathcal{S}=((X,d),(f_{i})_{i\in I})\text{\textit{.}}
\end{equation*}

\medskip

The \textit{fractal operator} $F_{\mathcal{S}}:P_{cp}(X)\rightarrow
P_{cp}(X) $, associated to the $\varphi $-$\max $-IFS $\mathcal{S}$, is
given by%
\begin{equation*}
F_{\mathcal{S}}(K)=\underset{i\in I}{\cup }f_{i}(K)
\end{equation*}%
for every $K\in P_{cp}(X)$.

\medskip

We say that the $\varphi $-$\max $-IFS $\mathcal{S}$ has attractor if $F_{%
\mathcal{S}}$ is a Picard operator. The fixed point of $F_{\mathcal{S}}$ is
called the \textit{attractor} of the system $\mathcal{S}$ and it is denoted
by $A_{\mathcal{S}}$.

\bigskip

\textbf{Iterated function systems consisting of convex contractions}

\bigskip

\textbf{Definition 2.13 }(see Definition 3.1 from [7])\textbf{.\ }\textit{An
iterated function system consisting of convex contractions (IFSCC for short)
is described by:}

-\textit{\ a complete metric space }$(X,d)$

-\textit{\ a finite family of continuous functions }$(f_{i})_{i\in I}$%
\textit{, where }$f_{i}:X\rightarrow X$,\textit{\ such that} \textit{for
every }$i,j\in I$\textit{\ there exist }$a_{ij},b_{ij},c_{ij}\in \lbrack
0,\infty )$\textit{\ satisfying the following two properties:}

$\qquad \alpha $) $a_{ij}+b_{ij}+c_{ij}\overset{def}{=}d_{ij}$ \textit{and} $%
\underset{i,j\in I}{\max }$ $d_{ij}<1$;

$\qquad \beta $) 
\begin{equation*}
d((f_{i}\circ f_{j})(x),(f_{i}\circ f_{j})(y))\leq
a_{ij}d(x,y)+b_{ij}d(f_{i}(x),f_{i}(y))+c_{ij}d(f_{j}(x),f_{j}(y))
\end{equation*}%
\textit{for every }$i,j\in I$\textit{\ and every }$x,y\in X$\textit{.}

\medskip

\textit{We denote such a system by }%
\begin{equation*}
\mathcal{S}=((X,d),(f_{i})_{i\in I})\text{\textit{.}}
\end{equation*}

\bigskip

\textbf{Remark 2.14.} \textit{Each IFSCC} \textit{is a }$\varphi $-$\max $%
\textit{-IFS.}

Indeed, just take $p=2$\ and the comparison function $\varphi :[0,\infty
)\rightarrow \lbrack 0,\infty )$\ given by 
\begin{equation*}
\varphi (t)=(\underset{i,j\in I}{\max }\text{ }d_{ij})t\text{,}
\end{equation*}%
\ for every $t\in \lbrack 0,\infty )$.

\bigskip

\textbf{Theorem 2.15 }(see Theorem 3.2 from [7])\textbf{.} \textit{Each IFSCC%
} $\mathcal{S}$ \textit{has attractor.}

\bigskip

\textbf{Theorem 2.16 }(see Theorem 3.2, ii) and Theorem 3.6 from [7])\textbf{%
.} \textit{Given an IFSCC} $\mathcal{S}=((X,d),(f_{i})_{i\in I})$,\textit{\
the function }$\pi :\Lambda (I)\rightarrow A_{\mathcal{S}}$\textit{\ -where }%
$A_{\mathcal{S}}$\textit{\ is the attractor of }$\mathcal{S}$\textit{-
defined by }%
\begin{equation*}
\pi (\omega )=a_{\omega }\text{,}
\end{equation*}%
\textit{for every\ }$\omega =\omega _{1}\omega _{2}...\omega _{n}\omega
_{n+1}...\in \Lambda (I)$\textit{, where }$a_{\omega }=\underset{%
n\rightarrow \infty }{\lim }f_{\omega _{1}...\omega _{n}}(x)$ \textit{for
every} $x\in X$\textit{, has the following properties:}

\textit{a)\ it is continuous;}

\textit{b) it\ is onto;}

\textit{c)} 
\begin{equation*}
\pi \circ \tau _{i}=f_{i}\circ \pi \text{,}
\end{equation*}%
\textit{for every }$i\in I$\textit{, where }$\tau _{i}:\Lambda
(I)\rightarrow \Lambda (I)$\textit{\ is given by }$\tau _{i}(\omega
)=i\omega $\textit{\ for every }$\omega \in \Lambda (I)$\textit{.}

\bigskip

The function $\pi $ described by the above theorem is called \textit{the
canonical projection from }$\Lambda (I)$ \textit{to} $A_{\mathcal{S}}$.

\bigskip

\textbf{The metric spaces} $(\mathcal{C}_{b},d_{u})$ \textbf{and} $(\mathcal{%
C},d_{u})$

\bigskip

Given a nonempty set $I$ and a metric space $(X,d)$, we consider the metric
spaces$(\mathcal{C}_{b},d_{u})$, where $\mathcal{C}_{b}=\{f:\Lambda
(I)\rightarrow X\mid f$ is continuous and bounded$\}$ and 
\begin{equation*}
d_{u}(f,g)=\underset{\omega \in \Lambda (I)}{\sup }d(f(\omega ),g(\omega ))
\end{equation*}%
for every $f,g\in \mathcal{C}_{b}$.

\bigskip

\textbf{Remark 2.17.}

\textit{a) The metric space }$(\mathcal{C}_{b},d_{u})$\textit{\ is complete
provided that }$(X,d)$ \textit{is complete.}

\textit{b) If }$I$\textit{\ is finite, then }$\mathcal{C}_{b}=\{f:\Lambda
(I)\rightarrow X\mid f$ \textit{is continuous}$\}$\textit{. In this case, we
denote }$(\mathcal{C}_{b},d_{u})$ \textit{by} $(\mathcal{C},d_{u})$\textit{.}

\newpage

\textbf{The operator} $G_{\mathcal{S}}:\mathcal{C}\rightarrow \mathcal{C}$ 
\textbf{associated to a }$\varphi $\textbf{-max-IFS} $\mathcal{S}$

\bigskip

Given a $\varphi $-$\max $-IFS $\mathcal{S}=((X,d),(f_{i})_{i\in I})$ and $%
g\in \mathcal{C}$, we can consider the function $G_{\mathcal{S},g}:\Lambda
(I)\rightarrow X$ described by the equality 
\begin{equation*}
G_{\mathcal{S},g}(\omega )=f_{\omega _{1}}(g(\omega _{2}...\omega _{n}\omega
_{n+1}...))\text{,}
\end{equation*}%
for every $\omega =\omega _{1}\omega _{2}...\omega _{n}\omega _{n+1}...\in
\Lambda (I)$.

\bigskip

\textbf{Lemma 2.18.} \textit{For every }$\varphi $\textit{-}$\max $\textit{%
-IFS} $\mathcal{S}$, \textit{the function} $G_{\mathcal{S},g}$ \textit{is
continuous.}

\textit{Proof}: Let us suppose that $\mathcal{S}=((X,d),(f_{i})_{i\in I})$.

\medskip

\textbf{Claim 1}.\textbf{\ }$\Lambda (I)=\underset{i\in I}{\cup }\tau
_{i}(\Lambda (I))$.

The justification of this claim is obvious.

\medskip

\textbf{Claim 2}. The set $\tau _{i}(\Lambda (I))$ is open for every $i\in I$%
.

\textit{Justification of claim 2}. We have 
\begin{equation*}
\{\theta \in \Lambda (I)\mid d_{\Lambda }(\theta ,\omega )<\frac{1}{2}%
\}\subseteq \tau _{i}(\Lambda (I))\text{,}
\end{equation*}
for every $\omega \in \tau _{i}(\Lambda (I))$ and every $i\in I$.

\medskip

\textbf{Claim 3}. The restriction of $G_{\mathcal{S},g}$ to $\tau
_{i}(\Lambda (I))$ is continuos for every $i\in I$.

\textit{Justification of claim 3}. For every $i\in I$, the restriction of $%
G_{\mathcal{S},g}$ to $\tau _{i}(\Lambda (I))$ is $f_{i}\circ g\circ R$,
where $R:\Lambda (I)\rightarrow \Lambda (I)$ is given by $R(\omega
_{1}\omega _{2}\omega _{3}...)=\omega _{2}\omega _{3}....$ for every $\omega
=\omega _{1}\omega _{2}\omega _{3}...\in \Lambda (I)$. Note that $R$ is
continuous since $d_{\Lambda }(R(\omega ),R(\theta ))=2d_{\Lambda }(\omega
,\theta )$ for all $\omega ,\theta \in \Lambda (I)$. As $f_{i}$ and $g$ are
continuous, the justification of this claim is done.

\medskip

Consequently, taking into account the above claims and Theorem 18.2, f),
from [9], we conclude that $G_{\mathcal{S},g}$ is continuous. $\square $

\bigskip

\textbf{Definition 2.19.} \textit{The operator }$G_{\mathcal{S}}:\mathcal{C}%
\rightarrow \mathcal{C}$\textit{\ associated to a }$\varphi $-$\max $\textit{%
-IFS }$\mathcal{S}$\textit{\ is given by }%
\begin{equation*}
G_{\mathcal{S}}(g)=G_{\mathcal{S},g}
\end{equation*}%
\textit{for each }$g\in \mathcal{C}$\textit{.}

\bigskip

\textbf{Remark 2.20. }\textit{For a }$\varphi $-$\max $-IFS\textit{\ }$%
\mathcal{S}=((X,d),(f_{i})_{i\in I})$\textit{, we have:}

\textit{a)}%
\begin{equation*}
G_{\mathcal{S}}^{[n]}(g)(\omega )=f_{\omega _{n}\omega _{n-1}...\omega
_{2}\omega _{1}}(g(\omega _{n+1}\omega _{n+2}...))\text{\textit{,}}
\end{equation*}%
\textit{for every }$n\in \mathbb{N}$\textit{\ and every }$\omega =\omega
_{1}\omega _{2}...\omega _{m}\omega _{m+1}...\in \Lambda (I)$.

\textit{b)} 
\begin{equation*}
G_{\mathcal{S}}^{[n]}(g)(\Lambda (I))=F_{\mathcal{S}}^{[n]}(g((\Lambda (I)))%
\text{,}
\end{equation*}%
\textit{for every }$n\in \mathbb{N}$.

\bigskip

\textbf{Proposition 2.21.} \textit{The operator }$G_{\mathcal{S}}:\mathcal{C}%
\rightarrow \mathcal{C}$\textit{\ associated to a }$\varphi $-$\max $\textit{%
-IFS }$\mathcal{S}$\textit{\ is continuous.}

\textit{Proof}. Let us suppose that $\mathcal{S=}(X,(f_{i})_{i\in I})$.

It suffices to prove that $\underset{n\rightarrow \infty }{\lim }d_{u}(G_{%
\mathcal{S}}(g_{n}),G_{\mathcal{S}}(g))=0$ for every $g_{n},g\in \mathcal{C}$
such that $\underset{n\rightarrow \infty }{\lim }d_{u}(g_{n},g)=0$.

We have%
\begin{equation}
g_{n}(\Lambda (I))\subseteq (\underset{n\in \mathbb{N}}{\cup }g_{n}(\Lambda
(I)))\cup g(\Lambda (I))\overset{not}{=}B_{1}\text{,}  \tag{1}
\end{equation}%
for every $n\in \mathbb{N}$.

Note that $B_{1}$ is compact.

Indeed, since $(\Lambda (I),d_{\Lambda })$ is compact (see Remark 2.7, c)
and $g_{n}\overset{u}{\rightarrow }g$, Remark 2.2, a), assures us that $%
\underset{n\rightarrow \infty }{\lim }H(g_{n}(\Lambda (I)),g(\Lambda (I)))=0$%
. Since $g_{n}(\Lambda (I))\in P_{cp}(X)$ and $g(\Lambda (I))\in P_{cp}(X)$
(because $(\Lambda (I),d_{\Lambda })$ is compact and the functions $g$ and $%
g_{n}$ are continuous), Remark 2.2, b), leads to the conclusion that $B_{1}$
is compact.

Now, let us consider $\varepsilon >0$.

As the continuos functions of the finite family $(f_{i})_{i\in I}$ are
uniformly continuous on the compact set $B_{1}$, there exists $\delta
_{\varepsilon }>0$\ such that 
\begin{equation}
d(f_{i}(x),f_{i}(y))<\varepsilon \text{,}  \tag{2}
\end{equation}%
provided that $i\in I$\ and $x,y\in B_{1}$, $d(x,y)<\delta _{\varepsilon }$.

Moreover, since $\underset{n\rightarrow \infty }{\lim }d_{u}(g_{n},g)=0$,
there exists $n_{\varepsilon }\in \mathbb{N}$ such that 
\begin{equation}
d_{u}(g_{n},g)<\delta _{\varepsilon }\text{,}  \tag{3}
\end{equation}%
for every $n\in \mathbb{N}$, $n\geq n_{\varepsilon }$. Therefore 
\begin{equation}
d(g_{n}(\omega )),g(\omega ))\leq d_{u}(g_{n},g)\overset{(3)}{<}\delta
_{\varepsilon }\text{,}  \tag{4}
\end{equation}%
for every $\omega \in \Lambda (I)$ and every $n\in \mathbb{N}$, $n\geq
n_{\varepsilon }$. Making use of $(1)$, we infer that $g_{n}(\omega
),g(\omega )\in B_{1}$, so 
\begin{equation}
d((f_{i}(g_{n}(\omega )),f_{i}(g(\omega )))\overset{(2)\&(4)}{<}\varepsilon 
\text{,}  \tag{5}
\end{equation}%
for all $\omega \in \Lambda (I)$, $i\in I$ and $n\in \mathbb{N}$, $n\geq
n_{\varepsilon }$. Finally we have%
\begin{equation*}
d_{u}(G_{\mathcal{S}}(g_{n}),G_{\mathcal{S}}(g))=\underset{\omega \in
\Lambda (I)}{\sup }d((G_{\mathcal{S}}(g_{n}))(\omega ),(G_{\mathcal{S}%
}(g))(\omega ))=
\end{equation*}%
\begin{equation*}
=\underset{i\in I,\omega \in \Lambda (I)}{\sup }d((G_{\mathcal{S}%
}(g_{n}))(i\omega ),(G_{\mathcal{S}}(g))(i\omega ))=\underset{i\in I,\omega
\in \Lambda (I)}{\sup }d((f_{i}(g_{n}(\omega )),f_{i}(g(\omega )))\overset{%
(5)}{\leq }\varepsilon \text{,}
\end{equation*}%
for every $n\in \mathbb{N}$, $n\geq n_{\varepsilon }$, i.e. $\underset{%
n\rightarrow \infty }{\lim }d_{u}(G_{\mathcal{S}}(g_{n}),G_{\mathcal{S}%
}(g))=0$. $\square $

\bigskip

\textbf{3.} \textbf{The main results}

\bigskip

\textbf{Theorem 3.1.} \textit{The operator} $G_{\mathcal{S}}$ \textit{%
associated to a }$\varphi $-$\max $\textit{-IFS} $\mathcal{S}$ \textit{is a
Picard operator.}

\textit{Proof}. Let us suppose that $\mathcal{S}=(X,(f_{i})_{i\in I})$.

We are going to prove that%
\begin{equation}
d_{u}(G_{\mathcal{S}}^{[p]}(g),G_{\mathcal{S}}^{[p]}(h))\leq \varphi (%
\underset{j\in \{0,1,...,p-1\}}{\max }d_{u}(G_{\mathcal{S}}^{[j]}(g),G_{%
\mathcal{S}}^{[j]}(h)))\text{,}  \tag{1}
\end{equation}%
for every $g,h\in \mathcal{C}$.

Indeed, we have%
\begin{equation*}
d_{u}(G_{\mathcal{S}}^{[p]}(g),G_{\mathcal{S}}^{[p]}(h))=\underset{\omega
\in \Lambda (I)}{\sup }d((G_{\mathcal{S}}^{[p]}(g))(\omega ),(G_{\mathcal{S}%
}^{[p]}(h))(\omega ))=
\end{equation*}%
\begin{equation*}
=\underset{\theta \in \Lambda _{p}(I),\omega \in \Lambda (I)}{\sup }d(((G_{%
\mathcal{S}}^{[p]}(g))(\theta \omega ),(G_{\mathcal{S}}^{[p]}(h))(\theta
\omega ))\overset{\text{Remark 2.20, a)}}{=}
\end{equation*}%
\begin{equation*}
=\underset{\theta \in \Lambda _{p}(I),\omega \in \Lambda (I)}{\sup }%
d(f_{\theta }(g(\omega )),f_{\theta }(h(\omega ))))\overset{\text{Definition
2.12}}{\leq }
\end{equation*}%
\begin{equation*}
\leq \underset{\omega \in \Lambda (I)}{\sup }\varphi (\underset{\sigma \in
V_{p}(I)}{\max }d(f_{\sigma }(g(\omega )),f_{\sigma }(h(\omega ))))\overset{%
\text{Definition 2.3}}{\leq }
\end{equation*}%
\begin{equation*}
\leq \varphi (\underset{\sigma \in V_{p}(I)}{\max }\underset{\omega \in
\Lambda (I)}{\sup }d(f_{\sigma }(g(\omega )),f_{\sigma }(h(\omega ))))=
\end{equation*}%
\begin{equation*}
=\varphi (\underset{j\in \{0,1,...,p-1\}}{\max }\underset{\sigma \in \Lambda
_{j}(I),\omega \in \Lambda (I)}{\sup }d(f_{\sigma }(g(\omega )),f_{\sigma
}(h(\omega ))))\overset{\text{Remark 2.20, a)}}{=}
\end{equation*}%
\begin{equation*}
=\varphi (\underset{j\in \{0,1,...,p-1\}}{\max }\underset{\omega \in \Lambda
(I)}{\sup }d(G_{\mathcal{S}}^{[j]}(g)(\omega ),G_{\mathcal{S}%
}^{[j]}(h)(\omega ))=
\end{equation*}%
\begin{equation*}
=\varphi (\underset{j\in \{0,1,...,p-1\}}{\max }d_{u}(G_{\mathcal{S}%
}^{[j]}(g),G_{\mathcal{S}}^{[j]}(h)))\text{,}
\end{equation*}%
for every $g,h\in \mathcal{C}$.

Taking into account Remark 2.17, a), Proposition 2.21 and $(1)$, based on
Theorem 2.6, we conclude that $G_{\mathcal{S}}$\textit{\ }is a Picard
operator. $\square $

\bigskip

\textbf{Theorem 3.2.} \textit{The fractal operator} $F_{\mathcal{S}}$\textit{%
\ associated to a} $\varphi $-$\max $\textit{-IFS} $\mathcal{S}$ \textit{is
a Picard operator.}

\textit{Proof}. Let us suppose that $\mathcal{S}=((X,d),(f_{i})_{i\in I})$.
According to Theorem 3.1, $G_{\mathcal{S}}$\ has a unique fixed point $%
g_{0}\in \mathcal{C}$.

\medskip

First we consider the case that $I$ has at least two elements.

\medskip

\textbf{Claim 1.} $g_{0}(\Lambda (I))\in P_{cp}(X)$\textit{.}

\textit{Justification of claim 1}. Since $(\Lambda (I),d_{\Lambda })$ is
compact (see Remark 2.7, c) and $g_{0}$ is continuous, we infer that $%
g_{0}(\Lambda (I))$ is a compact subset of $X$.

\medskip

\textbf{Claim 2.} $g_{0}(\Lambda (I))$ is a fixed point of\textit{\ }$F_{%
\mathcal{S}}$\textit{.}

\textit{Justification of claim 2}. Because $G_{\mathcal{S}}(g_{0})=g_{0}$,
we get $G_{\mathcal{S}}(g_{0})(i\omega )=g_{0}(i\omega )$, i.e. $%
f_{i}(g_{0}(\omega ))=g_{0}(\tau _{i}(\omega ))$ for every $i\in I$ and
every $\omega \in \Lambda (I)$. Consequently%
\begin{equation}
f_{i}\circ g_{0}=g_{0}\circ \tau _{i}\text{,}  \tag{1}
\end{equation}%
for every $i\in I$.

Then%
\begin{equation*}
g_{0}(\Lambda (I))=g_{0}(\underset{i\in I}{\cup }\tau _{i}(\Lambda (I)))=%
\underset{i\in I}{\cup }(g_{0}\circ \tau _{i})(\Lambda (I))\overset{(1)}{=}
\end{equation*}%
\begin{equation*}
=\underset{i\in I}{\cup }(f_{i}\circ g_{0})(\Lambda (I))=F_{\mathcal{S}%
}(g_{0}(\Lambda (I)))\text{,}
\end{equation*}%
i.e. $g_{0}(\Lambda (I))$ is a fixed point of\textit{\ }$F_{\mathcal{S}}$.

\medskip

\textbf{Claim 3.} $\underset{n\rightarrow \infty }{\lim }F_{\mathcal{S}%
}^{[n]}(K)=g_{0}(\Lambda (I))$ for every $K\in P_{cp}(X)$.

\textit{Justification of claim 3}. For every $g\in \mathcal{C}$, we have $%
\underset{n\rightarrow \infty }{\lim }G_{\mathcal{S}}^{[n]}(g)=g_{0}$, i.e. $%
G_{\mathcal{S}}^{[n]}(g)\overset{u}{\rightarrow }g_{0}$. Hence, since $%
(\Lambda (I),d_{\Lambda })$ is compact (see Remark 2.7, c), according to
Remark 2.2, a), we obtain $\underset{n\rightarrow \infty }{\lim }G_{\mathcal{%
S}}^{[n]}(g)(\Lambda (I))=g_{0}(\Lambda (I))$. In view of Remark 2.20, b),
we get 
\begin{equation}
\underset{n\rightarrow \infty }{\lim }F_{\mathcal{S}}^{[n]}(g((\Lambda
(I)))=g_{0}(\Lambda (I))\text{.}  \tag{2}
\end{equation}

Proposition 2.11 assures us that for every $K\in P_{cp}(X)$, there exists $%
g_{K}\in \mathcal{C}$ such that $g_{K}(\Lambda (I))=K$ and from $(2)$ we
infer that $\underset{n\rightarrow \infty }{\lim }F_{\mathcal{S}%
}^{[n]}(K)=g_{0}(\Lambda (I))$.

\medskip

\textbf{Claim 4.} $g_{0}(\Lambda (I))$ is the unique fixed point of $F_{%
\mathcal{S}}$.

\textit{Justification of claim 4}. If $A\in P_{cp}(X)$ is a fixed point of $%
F_{\mathcal{S}}$, then $F_{\mathcal{S}}^{[n]}(A)=A$ for every $n\in \mathbb{N%
}$, so $\underset{n\rightarrow \infty }{\lim }F_{\mathcal{S}}^{[n]}(A)=A$.
But, according to Claim 3, we have $\underset{n\rightarrow \infty }{\lim }F_{%
\mathcal{S}}^{[n]}(A)=g_{0}(\Lambda (I))$. Therefore, the uniqueness of the
limit of $(F_{\mathcal{S}}^{[n]}(A))_{n\in \mathbb{N}}$ implies that $%
A=g_{0}(\Lambda (I))$.

\medskip

Based on the above claims we conclude that $F_{\mathcal{S}}$ is a Picard
operator whose fixed point is $g_{0}(\Lambda (I))$.

\medskip

If $I$ has just one element, let us say $I=\{i\}$, then we consider $j\neq i$
and apply the above considerations to the system $\mathcal{S}^{^{\prime
}}=((X,d),\{f_{i},f_{j}\})$, where $f_{i}=f_{j}$. Hence $F_{\mathcal{S}%
^{^{\prime }}}$ is a Picard operator, so $F_{\mathcal{S}}$ is also a Picard
operator. $\square $

\bigskip

In view of Remark 2.14, we get the following:

\bigskip

\textbf{Theorem 3.3.} \textit{The operator} $G_{\mathcal{S}}$ \textit{%
associated to an IFSCC} $\mathcal{S}$ \textit{is a Picard operator.}

\bigskip

\textbf{Theorem 3.4.} \textit{The unique fixed point of the operator} $G_{%
\mathcal{S}}$ \textit{associated to an IFSCC} $\mathcal{S}%
=((X,d),(f_{i})_{i\in I})$ \textit{is} \textit{the canonical projection from 
}$\Lambda (I)$ \textit{to} $A_{\mathcal{S}}$, \textit{where }$A_{\mathcal{S}%
} $\textit{\ is the attractor of }$\mathcal{S}$.

\textit{Proof}. By virtue of Theorem 3.3, there exists a unique $%
g_{0}:\Lambda (I)\rightarrow X$ continuous such that $G_{\mathcal{S}%
}(g_{0})=g_{0}$. Consequently $G_{\mathcal{S}}(g_{0})(i\omega
)=g_{0}(i\omega )$, i.e.%
\begin{equation}
f_{i}(g_{0}(\omega ))=g_{0}(i\omega )\text{,}  \tag{1}
\end{equation}%
for every $i\in I$ and every $\omega \in \Lambda (I)$.

Therefore, using $(1)$, we get%
\begin{equation}
g_{0}(\omega _{1}...\omega _{n}\theta )=f_{\omega _{1}...\omega
_{n}}(g_{0}(\theta ))\text{,}  \tag{2}
\end{equation}%
for every $\theta $, $\omega =\omega _{1}\omega _{2}...\omega _{n}\omega
_{n+1}...\in \Lambda (I)$ and every $n\in \mathbb{N}$.

As $d_{\Lambda }(\omega _{1}...\omega _{n}\theta ,\omega )\leq \frac{1}{2^{n}%
}$ for every $n\in \mathbb{N}$, we infer that 
\begin{equation}
\underset{n\rightarrow \infty }{\lim }d_{\Lambda }(\omega _{1}...\omega
_{n}\theta ,\omega )=0\text{,}  \tag{3}
\end{equation}%
for every $\theta $, $\omega =\omega _{1}\omega _{2}...\omega _{n}\omega
_{n+1}...\in \Lambda (I)$.

On the one hand, from $(3)$, using the continuity of $g_{0}$, we obtain that 
$\underset{n\rightarrow \infty }{\lim }d(g_{0}(\omega _{1}...\omega
_{n}\theta ),g_{0}(\omega ))=0$, which, based on $(2)$, takes the following
form: 
\begin{equation}
\underset{n\rightarrow \infty }{\lim }d(f_{\omega _{1}...\omega
_{n}}(g_{0}(\theta )),g_{0}(\omega ))=0\text{,}  \tag{4}
\end{equation}%
for every $\theta $, $\omega =\omega _{1}\omega _{2}...\omega _{n}\omega
_{n+1}...\in \Lambda (I)$.

On the other hand, taking into account Theorem 2.16, we have%
\begin{equation}
\underset{n\rightarrow \infty }{\lim }d(f_{\omega _{1}...\omega
_{n}}(g_{0}(\theta )),\pi (\omega ))=0\text{,}  \tag{5}
\end{equation}%
for every $\theta $, $\omega =\omega _{1}\omega _{2}...\omega _{n}\omega
_{n+1}...\in \Lambda (I)$.

Relations $(4)$ and $(5)$ lead to the conclusion that $g_{0}=\pi $, i.e. the
fixed point of $G_{\mathcal{S}}$ is the canonical projection from $\Lambda
(I)$ to $A_{\mathcal{S}}$. $\square $

\bigskip

\textbf{4. Iterated function systems involving a possible infinite family of 
}$\varphi $\textbf{-max-contractions}

\bigskip

Theorem 3.1 is also valid for a class of iterated function systems involving
a possible infinite family of $\varphi $-$\max $-contractions. More
precisely, we have the following:

\bigskip

\textbf{Definition 4.1.\ }\textit{A possibly infinite iterated function
system consisting of }$\varphi $-$\max $\textit{-contractions (}$\varphi $-$%
\max $-\textit{PIIFS for short) is described by:}

-\textit{\ a complete metric space }$(X,d)$

-\textit{\ a family of continuous functions }$(f_{i})_{i\in I}$\textit{,
where }$f_{i}:X\rightarrow X$,\textit{\ having the following three
properties:}

\qquad a) $\underset{i\in I}{\cup }f_{i}(B)$\textit{\ is bounded, for every }%
$B\in P_{b,cl}(X)$\textit{;}

\qquad b) \textit{for each bounded and closed subset }$B$ \textit{of} $X$ 
\textit{and each }$\varepsilon >0$\textit{\ there exists }$\delta
_{B,\varepsilon }>0$\textit{\ such that }$d(f_{i}(x),f_{i}(y))<\varepsilon $%
\textit{\ provided that }$i\in I$\textit{\ and }$x,y\in B$\textit{, }$%
d(x,y)<\delta _{B,\varepsilon }$\textit{\ (i.e. the family }$(f_{i})_{i\in
I} $ \textit{is equal uniformly continuous on each bounded and closed subset
of} $X$\textit{).}

\qquad c)\textit{\ there exist a comparison function }$\varphi :[0,\infty
)\rightarrow \lbrack 0,\infty )$\textit{\ and }$p\in \mathbb{N}^{\ast }$ 
\textit{such that}%
\begin{equation*}
\underset{\omega \in \Lambda _{p}(I)}{\max }d(f_{\omega }(x),f_{\omega
}(y))\leq \varphi (\underset{\omega \in V_{p}(I)}{\max }d(f_{\omega
}(x),f_{\omega }(y)))\text{,}
\end{equation*}%
\textit{for every }$x,y\in X$\textit{.}

\textit{We denote such a system by }%
\begin{equation*}
\mathcal{S}=((X,d),(f_{i})_{i\in I})\text{\textit{.}}
\end{equation*}

\medskip

The \textit{fractal operator} $F_{\mathcal{S}}:P_{b,cl}(X)\rightarrow
P_{b,cl}(X)$, associated to the $\varphi $-$\max $-PIIFS $\mathcal{S}$, is
given by%
\begin{equation*}
F_{\mathcal{S}}(B)=\overline{\underset{i\in I}{\cup }f_{i}(B)}
\end{equation*}%
for every $B\in P_{cp}(X)$.

\bigskip

\textbf{Remark 4.2.} \textit{We can associate to a }$\varphi $\textit{-}$%
\max $\textit{-PIIFS }$\mathcal{S}$\textit{\ the operator }$G_{\mathcal{S}}:%
\mathcal{C}_{b}\rightarrow \mathcal{C}_{b}$\textit{\ described by }$G_{%
\mathcal{S}}(g)=G_{\mathcal{S},g}$.

Indeed, the continuity of\textit{\ }$G_{\mathcal{S},g}$ could be prove using
exactly the same arguments as those from the proof of Lemma 2.18 and the
justification of its boundedness is the following:

As$\ g$ is bounded, $(g\circ R)(\Lambda (I))$ is bounded and, in view of
property a) from Definition 4.1, $\underset{i\in I}{\cup }f_{i}((g\circ
R)(\Lambda (I)))=G_{\mathcal{S},g}(\Lambda (I)))$\textit{\ }is bounded, i.e. 
$G_{\mathcal{S},g}$ is bounded.

\bigskip

\textbf{Proposition 4.3.} \textit{The operator} $G_{\mathcal{S}}$ \textit{%
associated to a }$\varphi $-$\max $\textit{-PIIFS} $\mathcal{S}%
=(X,(f_{i})_{i\in I})$ \textit{is continuous.}

\textit{Proof}. The proof is similar with the one of Proposition 2.21,
except the justification of relations $(1)$ and $(2)$ from that proof which
should be replace with the following one:

Since $\underset{n\rightarrow \infty }{\lim }d_{u}(g_{n},g)=0$, there exists 
$n_{0}\in \mathbb{N}$ such that 
\begin{equation}
d_{u}(g_{n},g)\leq \frac{1}{2}\text{,}  \tag{*}
\end{equation}%
for every $n\in \mathbb{N}$, $n\geq n_{0}$.

We have%
\begin{equation}
g_{n}(\Lambda (I))\subseteq \overline{E_{1}(g(\Lambda (I)))}\overset{not}{=}%
B_{1}\text{,}  \tag{1}
\end{equation}%
for every $n\in \mathbb{N}$, $n\geq n_{0}$.

Indeed, for every $\omega \in \Lambda (I)$ and every $n\in \mathbb{N}$, $%
n\geq n_{0}$, we have $d(g_{n}(\omega ),g(\omega ))$\linebreak $\leq
d_{u}(g_{n},g)\overset{(\ast )}{\leq }\frac{1}{2}<1$, so $g_{n}(\omega )\in
\{y\in X\mid d(y,g(\omega ))<\frac{1}{2}\}$ and consequently $g_{n}(\omega
)\in \overline{E_{1}(g(\Lambda (I)))}$.

According to b) from the Definition 4.1, as $B_{1}$ is closed and bounded
(since $g(\Lambda (I))$ is bounded as $g$ is bounded), there exists $\delta
_{\varepsilon }>0$\ such that 
\begin{equation}
d(f_{i}(x),f_{i}(y))<\varepsilon \text{,}  \tag{2}
\end{equation}%
provided that $i\in I$\ and $x,y\in B_{1}$, $d(x,y)<\delta _{\varepsilon }$. 
$\square $

\bigskip

\textbf{Theorem 4.4.} \textit{The operator} $G_{\mathcal{S}}$ \textit{%
associated to a }$\varphi $-$\max $\textit{-PIIFS} $\mathcal{S}$ \textit{is
a Picard operator.}

\textit{Proof}. The proof is identical with the one of Theorem 3.1.

\bigskip

\textbf{Open problem}. \textit{Is it true that the fractal operator }$F_{%
\mathcal{S}}$\textit{\ associated to a }$\varphi $\textit{-}$\max $\textit{%
-PIIFS }$\mathcal{S}$\textit{\ is a Picard operator whose fixed point is the
image of the fixed point of }$G_{\mathcal{S}}$\textit{?}

\medskip

If the answer to this question is true, then the image of the fixed point of 
$G_{\mathcal{S}}$ would deserve to be called the attractor of the $\varphi $%
\textit{-}$\max $-PIIFS\textit{\ }$\mathcal{S}$.

\bigskip

\textbf{References}

\bigskip

[1] Barnsley, M.F.: Transformation between attractors of hyperbolic iterated
function systems (2007), arXiv:math/0703398v1.

[2] Brucks, K.M. and Bruin, H.: Topics from one-dimensional dynamics,
Cambridge University Press, 2004.

[3] Hille M.R.: Remarks on limit sets of infinite iterated function systems,
Monatsh. Math., \textbf{168} (2012), 215--237.

[4] Georgescu F.: IFSs consisting of generalized convex contractions, An. 
\c{S}tiin\c{t}. Univ. \textquotedblleft Ovidius\textquotedblright\ Constan%
\c{t}a, Ser. Mat., \textbf{25} (2017), 77-86.

[5] Mihail, A. and Miculescu, R.: The shift space for an infinite iterated
function system, Math. Rep., \textbf{61} (2009), 21-32.

[6] Miculescu, R. and Mihail, A.: A generalization of Matkowski's fixed
point theorem and Istr\u{a}\c{t}escu's fixed point theorem concerning convex
contractions, J. Fixed Point Theory Appl., 2017,
doi:10.1007/s11784-017-0411-7.

[7] Miculescu, R. and Mihail, A.: A generalization of Istr\u{a}\c{t}escu's
fixed point theorem for convex contractions, Fixed Point Theory, in print,
arXiv:1512.05490.

[8] Mihail, A.: The canonical projection between the shift space of an IIFS
and its attractor as a fixed point, Fixed Point Theory Appl., 2015, Paper
No. 75, 15 p.

[9] Munkres, J.R.: Topology, 2nd edition, Prentice Hall, Inc., Englewood
Cliffs, NJ, 2000.

[10] Willard, S.: General topology, Addison-Wesley Publishing Company, 1970.

\newpage

{\small Flavian Georgescu}

{\small Faculty of Mathematics and Computer Science}

{\small University of Pite\c{s}ti, Romania}

{\small T\^{a}rgul din Vale 1, 110040, Pite\c{s}ti, Arge\c{s}, Romania}

{\small E-mail: faviu@yahoo.com}

\bigskip

{\small Radu Miculescu}

{\small Faculty of Mathematics and Computer Science}

{\small Bucharest University, Romania}

{\small Str. Academiei 14, 010014, Bucharest}

{\small E-mail: miculesc@yahoo.com}

\bigskip

{\small Alexandru Mihail}

{\small Faculty of Mathematics and Computer Science}

{\small Bucharest University, Romania}

{\small Str. Academiei 14, 010014, Bucharest}

{\small E-mail: mihail\_alex@yahoo.com}

\end{document}